\documentclass[12pt]{article}%
\usepackage{amsfonts}
\usepackage{amssymb}
\usepackage{amsmath}%
\usepackage{graphicx}

\providecommand{\U}[1]{\protect\rule{.1in}{.1in}}

\newtheorem{theorem}{Theorem}

\newtheorem{corollary}[theorem]{Corollary}

\newtheorem{lemma}[theorem]{Lemma}

\newenvironment{proof}[1][Proof]{\noindent\textbf{#1.} }{{\hfill $\Box$ \\}}
\begin{document}

\title{Bases for free Lie superalgebras}
\author{Michael Vaughan-Lee}
\date{January 2024}
\maketitle

\begin{abstract}
We describe a basis for free Lie superalgebras which uses the theory of basic
commutators. The only description of bases for free Lie superalgbras that I
have found in the literature is in the book \emph{Infinite Dimensional Lie
Superalgebras} by Bahturin et al. [1]. Their bases make use of the theory of
Shirshov bases in free Lie algebras, and I believe that there is a case for
writing up an alternative approach using basic commutators. An additional reason
for publishing this note is that I use the basis described here in a
forthcoming paper where I prove that 5-Engel Lie algebras of characteristic
$p$ for $p>7$ are nilpotent of class at most 11.

\end{abstract}

\section{Introduction}

Let $L$ be a free Lie algebra over a commutative ring $R$ with 1, and let $L$
have free generating set $A$. Then $L$ is a free $R$-module, and there are
three standard ways described in the literature for obtaining a basis for $L$
over $R$. The first two are due to Lyndon \cite{lyndon54} and Shirshov
\cite{shirshov62}, and the third uses Philip Hall's collection process. To
obtain a Shirshov basis we proceed as follows. We let $A^{\ast}$ be the set of
all associative words $a_{1}a_{2}\ldots a_{n}$ ($n\geq1$) where $a_{1}%
,a_{2},\ldots,a_{n}\in A$ If $u=a_{1}a_{2}\ldots a_{n}\in A^{\ast}$ and
$v=b_{1}b_{2}\ldots b_{m}\in A^{\ast}$ then we set%
\[
uv=a_{1}a_{2}\ldots a_{n}b_{1}b_{2}\ldots b_{m},
\]
turning $A^{\ast}$ into a semigroup. We assume that there is a total order
\thinspace$<$ on the set $A$, and if $u,v\in A^{\ast}$ we let $u<v$ if $u$ is
lexicographically earlier than $v$. We say that $a_{1}a_{2}\ldots a_{n}$ is
\emph{regular} if%
\[
a_{1}a_{2}\ldots a_{n}>a_{i+1}a_{i+2}\ldots a_{n}a_{1}a_{2}\ldots a_{i}
\]
for all $i=1,2,\ldots,n-1$. (So elements of $A$ are regular.) We then define a
map $\pi$ from the set of regular words into $L$ as follows. If $a\in A$ then
we let $\pi(a)=a$. If $w$ is a regular word of length greater than 1 then we
write $w=uv$ where $v$ is chosen to have maximal length subject to $v$ being
regular and being a proper subword of $w$. It turns out that this choice of
$v$ implies that $u$ is also regular. Then we recursively define $\pi
(uv)=[\pi(u),\pi(v)]$. For a proof that the image of $\pi$ is a basis for $L$
(the Shirshov basis) see Chapter 2 of \cite{bahturin92}. The proof takes about
10 pages, and is extremely technical!

The Lyndon basis is defined similarly. A Lyndon word is an element $a_{1}%
a_{2}\ldots a_{n}\in A^{\ast}$ such that
\[
a_{1}a_{2}\ldots a_{n}<a_{i+1}a_{i+2}\ldots a_{n}a_{1}a_{2}\ldots a_{i}
\]
for all $i=1,2,\ldots,n-1$. (Elements of $A$ are Lyndon words.) We define a
map $\theta$ from the set of Lyndon words into $L$. If $a\in A$ then we let
$\theta(a)=a$. And if $w$ is a Lyndon word of length greater than 1 we write
$w=uv$ where $v$ is chosen to have maximal length subject to $v$ being a
Lyndon word, and $v$ being a proper subword of $w$. And as with the Shirshov
basis we recursively set $\theta(uv)=[\theta(u),\theta(v)]$. The set
$\{\theta(w)\,|\,w$ is a Lyndon word$\}$ is the Lyndon basis for $L$.

The third method of obtaining a basis for $L$ as a free $R$-module is to show
that the basic commutators (or basic Lie products) on the free generators $A$
form an $R$-module basis for $L$.

The basic commutators of weight one are the elements $a\in A$, which we assume
to be an ordered set.

The basic commutators of weight two are the elements $[a,b]$ where $a,b\in A$
and $a>b$. These are ordered arbitrarily among themselves, and so that they
follow basic commutators of weight one.

The basic commutators of weight three are the elements $[a,b,c]$ where
$a,b,c\in A$ and $a>b\leq c$. (We use the left-normed convention so that
$[a,b,c]$ denotes $[[a,b],c]$.) The basic commutators of weight three are
ordered arbitrarily among themselves, and so that they follow the basic
commutators of weight two.

In general, if $k>3$ the basic commutators of weight $k$ are the commutators
$[c,d]$ where for some $m,n$ such that $m+n=k$,

\begin{enumerate}
\item $c,d$ are basic commutators of weight $m,n$ respectively,

\item $c>d,$

\item if, in the definition of basic commutators of weight $m$, $c$ was
defined to be $[e,f]$ then $f\leq d$.
\end{enumerate}

The basic commutators of weight $k$ are then ordered arbitrarily among
themselves, and so that they follow the basic commutators of weight $k-1$.

A proof that the basic commutators form an $R$-module basis for $L$ can be
found in my book \cite{vlee93b}.

\section{Lie superalgebras}

A Lie superalgebra is a $\mathbb{Z}_{2}$-graded algebra $L=L_{0}\oplus L_{1}$
with a bilinear product $[,]$ such that
\begin{align*}
\lbrack L_{0},L_{0}],\,[L_{1},L_{1}]  & \leq L_{0},\\
\lbrack L_{0},L_{1}],\,[L_{1},L_{0}]  & \leq L_{1}.
\end{align*}
Elements in $L_{0}$ are said to be even elements, and elements in $L_{1}$ are
said to be odd. If $a$ is even then we set $|a|=0$, and if $a$ is odd then we
set $|a|=1$. Odd elements and even elements are said to be homogeneous.
Finally, the product $[,]$ must satisfy the following relations for all
homogeneous elements $a,b,c$.%
\[
\lbrack b,a]=-(-1)^{|a|.|b|}[a,b].
\]%
\[
(-1)^{|a|.|c|}[a,[b,c]]+(-1)^{|b|.|a|}[b,[c,a]]+(-1)^{|c|.|b|}[c,[a,b]]=0.
\]
It is helpful to note that these relations imply that if $a,b,c$ are
homogeneous elements then%
\[
\lbrack a,[b,c]]=[a,b,c]-(-1)^{|b|.|c|}[a,c,b].
\]
We also add in the requirement that $[a,a]=0$ for even elements, and the
requirement that $[a,a,a]=0$ for odd elements. (These extra requirements are
redundant if 2 and 3 are invertible.)

Now let $L=L_{0}\oplus L_{1}$ be a free Lie superalgebra over a ring $R$ with
1, freely generated by an ordered set $(A,<)$. We assume that the elements of
$A$ are homogeneous, some even and some odd. Let $S=\{\pi(w)\,|\,w$ is a
regular word in $A^{\ast}\}$, and let $(C,<)$ be a (complete) set of basic
commutators on $A$. Note that the definitions of $S$ and $C$ are independent
of the $\mathbb{Z}_{2}$-grading on $L$. Note also that the elements of $S$ and
$C$ are homogeneous elements of $L$.

\begin{theorem}
[{Bahturin et al. [1]}]The free Lie superalgebra $L$ is a free $R$-module with
basis $S\cup\{[w,w]\,|\,w\in S$ is odd$\}$.
\end{theorem}

\begin{theorem}
The free Lie superalgebra $L$ is a free $R$-module with basis 
\[
C\cup\{[w,w]\,|\,w\in \text{C is odd}\}.
\]
\end{theorem}

\section{Proof of Theorem 2}

Let $L=L_{0}\oplus L_{1}$ be a Lie superalgebra over a commutative ring $R$
with 1, and suppose that $L$ is generated by a set $A$ whose elements are all
homogeneous. (We are not assuming here that $L$ is free.) Let $(C,<)$ be a
complete set of basic commutators on the generators in $A$. As a preliminary
step to proving Theorem 2 we show that $C\cup\{[c,c]\,|\,c\in C$ is odd$\}$
spans $L$ as an $R$-module.

\begin{lemma}
Let $a_{1},a_{2},\ldots,a_{n}$ be homogeneous elements in $L$ and let $c$ be a
Lie product of $a_{1},a_{2},\ldots,a_{n}$ in some order with some bracketing.
Then $c$ is a $\mathbb{Z}$-linear combination of left-normed Lie products
$[a_{1},a_{2\sigma},a_{3\sigma},\ldots,a_{n\sigma}]$ where $\sigma$ is a
permutation of $\{2,3,\ldots,n\}$.
\end{lemma}

\begin{proof}
The proof is by induction on $n$, the case $n=1$ being trivial, and the case
$n=2$ following from the fact that $[a_{2},a_{1}]=\pm\lbrack a_{1},a_{2}]$. So
suppose that $n>2$, and suppose that our claim holds true for smaller values
of $n$. Let $a_{1},a_{2},\ldots,a_{n}$ be homogeneous elements in $L$ and let
$c$ be a Lie product of $a_{1},a_{2},\ldots,a_{n}$ in some order with some
bracketing. Let $c=[u,v]$. Since $[u,v]=\pm\lbrack v,u]$ it is sufficient to
consider the case when $u$ involves $a_{1}$. We next show that it is
sufficient to consider the case when $v$ has weight one as a Lie product in
(some of) $a_{2},a_{3},\ldots,a_{n}$. If $v$ has weight greater than 1 we
write $v=[x,y]$ and then%
\[
\lbrack u,v]=[u,[x,y]]=[[u,x],y]\pm\lbrack\lbrack u,y],x].
\]
So $[u,v]$ is a $\mathbb{Z}$-linear combination of two Lie products
$[u^{\prime},v^{\prime}]$ where $u^{\prime}=[u,x]$ or $[u,y]$, and where the
weight of $v^{\prime}$ is smaller than the weight of $v$. Repeating this
argument as necessary we see that $c=[u,v]$ is a $\mathbb{Z}$-linear
combination of Lie products $[u^{\prime},v^{\prime}]$ where $u^{\prime}$
involves $a_{1}$ and where $v^{\prime}$ has weight one. By induction
$u^{\prime}$ is a $\mathbb{Z}$-linear combination of left-normed Lie products
with first entry $a_{1}$, and this completes the proof of the lemma.
\end{proof}

Now let $c$ be a Lie product of weight $n$ in elements $a_{1},a_{2}%
,\ldots,a_{n}$ from the generating set $A$, and suppose that $a_{1}\leq a_{i}$
for $i=2,3,\ldots,n$. By Lemma 3 we can express $c$ as a $\mathbb{Z}$-linear
combination of left-normed Lie products $[c_{1},c_{2},\ldots,c_{n}]$ where
$c_{1}=a_{1}$ and $c_{2},c_{3},\ldots,c_{n}$ is a permutation of $a_{2}%
,a_{3},\ldots,a_{n}$.

Suppose first that $c_{1}<c_{2}$. Then $[c_{2},c_{1}]$ is a basic commutator
in $C$, and
\[
\lbrack c_{1},c_{2},\ldots,c_{n}]=\pm\lbrack\lbrack c_{2},c_{1}],c_{3}%
,\ldots,c_{n}].
\]
If we set $d_{2}=[c_{2},c_{1}]$, $d_{i}=c_{i}$ for $i=3,\ldots,n$, then we see
that
\[
\lbrack c_{1},c_{2},\ldots,c_{n}]=\pm\lbrack d_{2},d_{3},\ldots,d_{n}]
\]
and that $[d_{2},d_{3},\ldots,d_{n}]$ is a Lie product of basic commutators
$d_{2},d_{3},\ldots,d_{n}$ which have the property that if $2\leq i\leq n$ and
if $d_{i}=[e,f]$ in its definition as a basic commutator, then $f\leq d_{j}$
for $j=2,3,\ldots,n$.

Next suppose that $c_{1}=c_{2}$. If $c_{1}$ is even then $[c_{1},c_{2}]=0$, so
we may assume that $c_{1}$ is odd. If $n=2$ then $[c_{1},c_{2},\ldots
,c_{n}]=[c_{1},c_{1}]$ where $c_{1}$ is an odd basic commutator. If $n\geq3$
then we may suppose that $c_{3}>c_{1}$ since $[c_{1},c_{1},c_{1}]=0$. But then%
\[
\lbrack c_{1},c_{2},c_{3}]=[[c_{1},c_{1}],c_{3}]=-[c_{3},[c_{1},c_{1}%
]]=-2[c_{3},c_{1},c_{1}],
\]
and $[c_{3},c_{1},c_{1}]$ is a basic commutator in $C$. In this case, if we
set $d_{3}=[c_{3},c_{1},c_{1}]$, $d_{i}=c_{i}$ for $i=4,\ldots,n$, then we see
that
\[
\lbrack c_{1},c_{2},\ldots,c_{n}]=-2[d_{3},d_{4},\ldots,d_{n}]
\]
and that $[d_{3},d_{4},\ldots,d_{n}]$ is a Lie product of basic commutators
$d_{3},d_{4},\ldots,d_{n}$ which have the property that if $3\leq i\leq n$ and
if $d_{i}=[e,f]$ in its definition as a basic commutator, then $f\leq d_{j}$
for $j=3,4,\ldots,n$.

So we see that if $c$ is a Lie product of weight $n$ in elements $a_{1}%
,a_{2},\ldots,a_{n}$ from the generating set $A$, then $c$ is a $\mathbb{Z}%
$-linear combination of elements of the following three types:

\begin{enumerate}
\item a basic commutator in $C$,

\item $[c,c]$ where $c$ is an odd basic commutator in $C$,

\item $[c_{1},c_{2},\ldots,c_{k}]$ ($1<k<n)$ where $c_{1},c_{2},\ldots,c_{k}$
are basic commutators in $C$ with the property that if $1\leq i\leq k $ and if
$c_{i}$ is defined to be $[e,f]$ then $f\leq c_{j}$ for $j=1,2,\ldots,k$.
\end{enumerate}

Consider a Lie product of type 3. By Lemma 3 we can write it as a $\mathbb{Z}
$-linear combination of Lie products $[d_{1},d_{2},\ldots,d_{k}]$ where
$d_{1},d_{2},\ldots,d_{k}$ is a permutation of $c_{1},c_{2},\ldots,c_{k}$ and
where $d_{1}\leq d_{i}$ for $i=1,2,\ldots,k$. Let $[d_{1},d_{2},\ldots,d_{k}]$
be one of the Lie products in this linear combination.

If $d_{1}<d_{2}$ then $[d_{2},d_{1}]$ is a basic commutator in $C$ and
\[
\lbrack d_{1},d_{2},\ldots,d_{k}]=\pm\lbrack\lbrack d_{2},d_{1}],d_{3}%
,\ldots,d_{k}].
\]
Setting $e_{2}=[d_{2},d_{1}]$, $e_{i}=d_{i}$ for $i=3,4,\ldots,k$, we see that
$[[d_{2},d_{1}],d_{3},\ldots,d_{k}]$ is a left-normed Lie product of length
$k-1$ in basic commutators $e_{2},e_{3},\ldots,e_{k}$ with the property that
if $2\leq i\leq k$ and if $e_{i}$ is defined to be $[e,f]$ then $f\leq e_{j}$
for $j=2,3,\ldots,k$.

If $d_{1}=d_{2}$ then (as above) we see that either $[d_{1},d_{2},\ldots
,d_{k}]=0$, or if $k=2$ then $[d_{1},d_{2},\ldots,d_{k}]=[c,c]$ where $c$ is
an odd basic commutator in $C$, or if $k>2$ then $d_{3}>d_{1}$, $[d_{3}%
,d_{1},d_{1}]$ is a basic commutator in $C$, and%
\[
\lbrack d_{1},d_{2},\ldots,d_{k}]=-2[[d_{3},d_{1},d_{1}],d_{4},\ldots,d_{k}].
\]
If we set $e_{3}=[d_{3},d_{1},d_{1}]$, $e_{i}=d_{i}$ for $i=4,\ldots,k$, then
$[[d_{3},d_{1},d_{1}],d_{4},\ldots,d_{k}]$ is a left-normed Lie product of
length $k-2$ in basic commutators $e_{3},e_{4},\ldots,e_{k}$ with the property
that if $3\leq i\leq k$ and if $e_{i}$ is defined to be $[e,f]$ then $f\leq
e_{j}$ for $j=3,\ldots,k$.

Continuing in this way by reverse induction on $k$ we see that any Lie product
of elements from the generating set $A$ is a $\mathbb{Z}$-linear combination
of elements of type 1 or type 2 above. This proves the following lemma.

\begin{lemma}
Let $L=L_{0}\oplus L_{1}$ be a Lie superalgebra over a commutative ring $R$
with $1$, and suppose that $L$ is generated by a set $A$ whose elements are all
homogeneous. Let $(C,<)$ be a complete set of basic commutators on the
generators in $A$. Then
\[
C\cup\{[c,c]\,|\,c\in C \text{ is odd}\}
\]
spans $L$ as an $R$-module.
\end{lemma}

To complete the proof of Theorem 2 we need to show that in the case when $L$
is freely generated by $A$ then $C\cup\{[c,c]\,|\,c\in C$ is odd$\}$ is
linearly independent over $R$. It is sufficient to prove that this is the case
when $L$ is a free Lie superalgbra over the integers $\mathbb{Z}$, as this
will imply that $L\otimes_{\mathbb{Z}}R$ is a free Lie superalgebra over $R$.
We let $A=A_{0}\oplus A_{1}$ be the free associative superalgebra over
$\mathbb{Z}$, freely generated by $a_{1},a_{2},\ldots,a_{r}$, where these free
generators are homogeneous (some even and some odd). So if $a$ and $b$ are two
elements in $A$ and either both lie in $A_{0}$ or both lie in $A_{1}$ then
$ab\in A_{0}$, and if $a$ and $b$ are two elements in $A$ with one element in
$A_{0}$ and the other in $A_{1}$ the $ab\in A_{1}$. We define a bilinear bracket
product $[,]$ on $A$ by setting%
\[
\lbrack a,b]=ab-(-1)^{|a|.|b|}ba
\]
for homogeneous elements $a,b\in A$. It is straightforward to see that this
turns $A$ into a Lie superalgebra. We show that if $L$ is the Lie
sub-superalgebra of $A$ generated by $a_{1},a_{2},\ldots,a_{r}$ then $L$ is
free, and has a $\mathbb{Z}$-module basis $C\cup\{[c,c]\,|\,c\in C$ is odd$\}$
where $C$ is a complete set of basic commutators on the generators.

Let $A_{n}$ be the $\mathbb{Z}$-module spanned by all products $a_{i_{1}%
}a_{i_{2}}\ldots a_{i_{n}}$ of length $n$ in the generators $a_{1}%
,a_{2},\ldots,a_{r}$. Clearly $A_{n}$ is a free $\mathbb{Z}$-module of rank
$r^{n}$, with these products forming a basis. We define a \emph{basic product}
in $A$ to be an element of the form $c_{1}c_{2}\ldots c_{m}$ ($m\geq1$) where
$c_{1},c_{2},\ldots,c_{m}$ are basic commutators in $C$, and where $c_{1}\leq
c_{2}\leq\ldots\leq c_{m}$. And we define the weight of a basic product
$c_{1}c_{2}\ldots c_{m}$ to be $\sum_{i=1}^{m}$wt\thinspace($c_{i}$).

\begin{lemma}
There are $r^{n}$ basic products of weight $n$, and they span $A_{n}$.
\end{lemma}

\begin{proof}
We use the Hall collection process to express the generators $a_{i_{1}%
}a_{i_{2}}\ldots a_{i_{n}}$ of $A_{n}$ as linear combinations of basic
products. Consider a generator $a_{i_{1}}a_{i_{2}}\ldots a_{i_{n}}$. First we
collect entries $a_{1}$ in the product towards the left. We look for a subword
$a_{i_{j}}a_{i_{j+1}}$ where $i_{j}>1$ and $i_{j+1}=1$. If there are no such
subwords then all entries $a_{1}$ are collected to the left, and%
\[
a_{i_{1}}a_{i_{2}}\ldots a_{i_{n}}=a_{1}^{m}a_{j_{1}}a_{j_{2}}\ldots
a_{j_{n-m}}
\]
for some $m\geq0$, and some $j_{1},j_{2},\ldots,j_{n-m}>1$. If there is a
subword $a_{i_{j}}a_{i_{j+1}}$ of this form then $a_{i_{j}}a_{i_{j+1}%
}=[a_{i_{j}},a_{1}]+\varepsilon a_{1}a_{i_{j}}$ where $\varepsilon=\pm1$, and
$[a_{i_{j}},a_{1}]\in C$. We then replace $a_{i_{1}}a_{i_{2}}\ldots a_{i_{n}}$
by%
\[
a_{i_{1}}a_{i_{2}}\ldots a_{i_{j-1}}[a_{i_{j}},a_{1}]a_{i_{j+2}}\ldots
a_{i_{n}}+\varepsilon a_{i_{1}}a_{i_{2}}\ldots a_{i_{j-1}}a_{1}a_{i_{j}%
}a_{i_{j+2}}\ldots a_{i_{n}}.
\]
We continue collecting $a_{1}$ to the left in each of the products%
\[
a_{i_{1}}a_{i_{2}}\ldots a_{i_{j-1}}[a_{i_{j}},a_{1}]a_{i_{j+2}}\ldots
a_{i_{n}},\;a_{i_{1}}a_{i_{2}}\ldots a_{i_{j-1}}a_{1}a_{i_{j}}a_{i_{j+2}%
}\ldots a_{i_{n}}
\]
in turn. If $i_{j+2}=1$ then
\[
\lbrack a_{i_{j}},a_{1}]a_{i_{j+2}}=[a_{i_{j}},a_{1}]a_{1}=[a_{i_{j}}%
,a_{1},a_{1}]+\eta a_{1}[a_{i_{j}},a_{1}]
\]
where $\eta=\pm1$, and we substitute $[a_{i_{j}},a_{1},a_{1}]+\eta
a_{1}[a_{i_{j}},a_{1}]$ for $[a_{i_{j}},a_{1}]a_{i_{j+2}}$ in the product%
\[
a_{i_{1}}a_{i_{2}}\ldots a_{i_{j-1}}[a_{i_{j}},a_{1}]a_{i_{j+2}}\ldots
a_{i_{n}}.
\]
If $i_{j+2}>1$ then we search $a_{i_{1}}a_{i_{2}}\ldots a_{i_{j-1}}[a_{i_{j}%
},a_{1}]a_{i_{j+2}}\ldots a_{i_{n}}$ for a subword $a_{i}a_{1}$ with $i>1$,
and if we find one then we substitute $[a_{i},a_{1}]\pm a_{1}a_{i}$ for
$a_{i}a_{1}$, as before. We apply the same procedure to
\[
a_{i_{1}}a_{i_{2}}\ldots a_{i_{j-1}}a_{1}a_{i_{j}}a_{i_{j+2}}\ldots a_{i_{n}%
}.
\]

We introduce the notation $[b,_{r}a]$ $(r=0,1,2,\ldots)$ for repeated
commutators, setting $[b,_{0}a]=b$, $[b,_{1}a]=[b,a]$, $[b,_{2}a]=[b,a,a]\,$,
and so on. At any stage in the process of collecting entries $a_{1}$ to the
left we have an expression for $a_{i_{1}}a_{i_{2}}\ldots a_{i_{n}}$ as a
linear combination of products of the form%
\begin{equation}
a_{1}^{e_{1}}[a_{j_{1}},_{m_{1}}a_{1}]a_{1}^{e_{2}}[a_{j_{2}},_{m_{2}}%
a_{1}]a_{1}^{e_{3}}\ldots\lbrack a_{j_{k}},_{m_{k}}a_{1}]a_{1}^{e_{k+1}}%
\end{equation}
where (for some $k$) $e_{1},e_{2},\ldots,e_{k+1},m_{1},m_{2},\cdots,m_{k}%
\geq0$, and where \thinspace$j_{1},j_{2},\ldots,j_{k}>1$. (Note that the
commutators $[a_{j},_{m}a_{1}]$ are basic commutators.) If $e_{2}%
=\ldots=e_{k+1}=0$ then all entries $a_{1}$ in (1) are collected to the left.
If $e_{s+1}>0$ then
\[
\lbrack a_{j_{s}},_{m_{s}}a_{1}]a_{1}^{e_{s+1}}=[a_{j_{s}},_{m_{s}}a_{1}%
]a_{1}a_{1}^{e_{s+1}-1}=[a_{j_{s}},_{m_{s}+1}a_{1}]a_{1}^{e_{s+1}%
-1}+\varepsilon a_{1}[a_{j_{s}},_{m_{s}}a_{1}]a_{1}^{e_{s+1}-1}%
\]
for some $\varepsilon=\pm1$, and we substitute 
\[[a_{j_{s}},_{m_{s}+1}a_{1}]a_{1}^{e_{s+1}-1}+\varepsilon %
a_{1}[a_{j_{s}},_{m_{s}}a_{1}]a_{1}^{e_{s+1}-1}
\]
for $[a_{j_{s}},_{m_{s}}a_{1}]a_{1}^{e_{s+1}}$ in (1).
Continuing in this way we eventually obtain an expression for $a_{i_{1}%
}a_{i_{2}}\ldots a_{i_{n}}$ as a linear combination of products
\begin{equation}
a_{1}^{e_{1}}[a_{j_{1}},_{m_{1}}a_{1}][a_{j_{2}},_{m_{2}}a_{1}]\ldots\lbrack
a_{j_{k}},_{m_{k}}a_{1}]
\end{equation}
where all entries $a_{1}$ are collected to the left. Now the spanning elements
$a_{i_{1}}a_{i_{2}}\ldots a_{i_{n}}$ for $A_{n}$ can all be expressed in the
form%
\begin{equation}
a_{1}^{e_{1}}a_{j_{1}}a_{1}^{m_{1}}a_{j_{2}}a_{1}^{m_{2}}\ldots a_{j_{k}}%
a_{1}^{m_{k}}%
\end{equation}
where $e_{1},m_{1},\ldots,m_{k}\geq0$, $j_{1},j_{2},\ldots,j_{k}>1$ and there
is a natural 1-1 correspondence between the set of all possible products of
weight $n$ of the form (2) and the set of all possible products of weight $n$
of the form (3). So there are $r^{n}$ possible products of weight $n$ of the
form (2), where all entries $a_{1}$ have been collected to the left.

We illustrate this by applying the Hall collection to a product $baba$ where
$a$ and $b$ are even generators of $A$ and $b>a$.%
\[
baba=[b,a]ba+abba.
\]
And%
\begin{align*}
\lbrack b,a]ba  & =[b,a][b,a]+[b,a]ab=[b,a][b,a]+[b,a,a]b+a[b,a]b,\\
abba  & =ab[b,a]+abab=ab[b,a]+a[b,a]b+aabb.
\end{align*}
So%
\[
baba=aabb+ab[b,a]+2a[b,a]b+[b,a][b,a]+[b,a,a]b.
\]
Note that there are six possible products in $a$ and $b$ which have weight two
in each of $a$ and $b$:%
\[
aabb,\,abab,\,abba,\,baab,\,baba,\,bbaa.
\]
And there are six possible products of this weight in which all entries $a$
have been collected to the left:%
\[
aabb,\,a[b,a]b,\,ab[b,a],\,[b,a,a]b,[\,b,a][b,a],\,b[b,a,a].
\]

Once we have an expression for $a_{i_{1}}a_{i_{2}}\ldots a_{i_{n}}$ as a
linear combination of products where all the entries $a_{1}$ have been
collected to the left, we collect entries $a_{2}$ in these products to the
left, and then entries $a_{3}$, and so on. There are only finitely many basic
commutators in $C$ which have weight at most $n$. Let these be $b_{1}%
,b_{2},\ldots,b_{m}$ with $b_{1}<b_{2}<\ldots<b_{m}$. Suppose we have reached
the point where all entries $b_{1},b_{2},\ldots,b_{k}$ have been collected to
the left, but no $b_{i}$ with $i>k$ has been collected to the left. Then we
have an expression for $a_{i_{1}}a_{i_{2}}\ldots a_{i_{n}}$ as a linear
combination of products of the form%
\begin{equation}
b_{1}^{e_{1}}b_{2}^{e_{2}}\ldots b_{k}^{e_{k}}b_{k+1}^{e_{k+1}}c_{1}%
b_{k+1}^{r_{1}}c_{2}b_{k+1}^{r_{2}}\ldots c_{t}b_{k+1}^{r_{t}}%
\end{equation}
where $e_{1},e_{2},\ldots,e_{k+1},r_{1},r_{2},\ldots,r_{t}\geq0$, and where
$c_{1},c_{2},\ldots,c_{t}$ are basic commutators in $C$ with $c_{i}>b_{k+1}$
for $i=1,2,\ldots,t$, and where if $c_{i}$ is defined to be $[e,f]$ then
$f\leq b_{k}$. We assume by induction that there are $r^{n}$ possible products
of the form (4) in $A_{n}$. We then collect entries $b_{k+1}$ to the left,
introducing basic commutators $[c_{i},_{r}b_{k+1}]$ $(r>0)$. When all entries
$b_{k+1}$ have been collected to the left then we will have an expression for
$a_{i_{1}}a_{i_{2}}\ldots a_{i_{n}}$ as a linear combination of products of
the form%
\begin{equation}
b_{1}^{e_{1}}b_{2}^{e_{2}}\ldots b_{k}^{e_{k}}b_{k+1}^{e_{k+1}}[c_{1},_{r_{1}%
}b_{k+1}][c_{2},_{r_{2}}b_{k+1}]\ldots\lbrack c_{t},_{r_{t}}b_{k+1}].
\end{equation}
There is a 1-1 correspondence between the set of all possible products of
weight $n$ of the form (4) and the set of all possible products of weight $n$
of the form (5). So we may assume that there are $r^{n}$ possible products of
weight $n$ of the form (5). Continuing in this way we eventually obtain an
expression for $a_{i_{1}}a_{i_{2}}\ldots a_{i_{n}}$ as a linear combination of
products of weight $n$ in which all entries $b_{1},b_{2},\ldots,b_{m}$ have
been collected to the left. These must all be basic products, and by induction
the number of possible basic products of weight $n$ is $r^{n}$.
\end{proof}

\begin{corollary}
The set $C\cup\{[c,c]\,|\,c\in C$ is odd$\}$ is linearly independent.
\end{corollary}

\begin{proof}
Since there are $r^{n}$ basic products of of weight $n$, and since they span
$A_{n}$, they must form a $\mathbb{Z}$-module basis for $A_{n}$ and must be
linearly independent. If $c\in C$ has weight $n$ then $c$ is a basic product
of weight $n$. And if $c\in C$ is odd and $[c,c]$ has weight $n$ then
$[c,c]=2c^{2}$ where $c^{2}$ is a basic product of weight $n$. So 
\[
\{c\in C\,|\,\text{wt}(c)=n\}\cup\{[c,c]\,|\,c\text{ is odd, wt}([c,c])=n\}
\]
is linearly independent. Since
\[
A=A_{1}\oplus A_{2}\oplus\ldots\oplus A_{n}\oplus\ldots
\]
this implies that $C\cup\{[c,c]\,|\,c\in C$ is odd$\}$ is linearly independent.
\end{proof}

Now let $M~$\ be a free Lie superalgebra over $\mathbb{Z}$, freely generated
by $x_{1},x_{2},\ldots,x_{r}$, and let $A=A_{0}\oplus A_{1}$ be the free
associative superalgebra over $\mathbb{Z}$, freely generated by $a_{1}%
,a_{2},\ldots,a_{r}$, where $|a_{i}|=|x_{i}|$ for $i=1,2,\ldots,r$. As above
we define a bracket product $[,]$ on $A$ by setting%
\[
\lbrack a,b]=ab-(-1)^{|a|.|b|}ba
\]
for homogeneous elements $a,b\in A$. And as above, we let $L$ be the Lie
sub-superalgebra of $A$ generated by $a_{1},a_{2},\ldots,a_{r}$. Let $D$ be a
complete set of basic commutators on $x_{1},x_{2},\ldots,x_{r}$, where we
assume that $x_{1}<x_{2}<\ldots<x_{r}$. By Lemma 4, $M$ is spanned by
\[
D\cup\{[d,d]\,|\,d\in D\text{ is odd}\}.
\]
There is a homomorphism from $M$ onto
$L$ mapping $x_{i}$ to $a_{i}$ for $i=1,2,\ldots,r$, and this homomorphism
maps $D$ onto a complete set of basic commutators on $a_{1},a_{2},\ldots
,a_{r}$. The image of $D\cup\{[d,d]\,|\,d\in D$ is odd$\}$ under this
homomorphism is linearly independent, and so $D\cup\{[d,d]\,|\,d\in D$ is
odd$\}$ is linearly independent. This proves Theorem 2 for free Lie
superalgebras of finite rank. But if Theorem 2 holds for free Lie superalgebras of
finite rank then it certainly holds true for free Lie superalgebras of arbitrary rank.

\end{document}